\documentclass[11pt,reqno]{article}
\usepackage[all]{xy}
\usepackage{amsmath}
\usepackage{amsthm}
\usepackage{amscd}
\usepackage{amssymb}

\input{epsf.tex}

\numberwithin{equation}{section}

\newcommand{\nb}[1]{#1\nobreakdash-}

\theoremstyle{definition}

\theoremstyle{plain}
\newtheorem{theorem}{Theorem}
\newtheorem{proposition}[theorem]{Proposition}
\newtheorem{corollary}[theorem]{Corollary}
\newtheorem{problem}[theorem]{Problem}

\newcounter{remarks}
{\paragraph*{Remarks}\smallskip
     \begin{list}{\arabic{remarks}. }{\usecounter{remarks}%
          \setlength{\leftmargin}{0in}%
          \setlength{\rightmargin}{0in}%
          \setlength{\labelsep}{0pt}%
          \setlength{\labelwidth}{0pt}%
          \setlength{\listparindent}{0pt}%
     }
}
{
\end{list}
}

\newcommand\inv\inverse

\DeclareMathOperator{\Ends}{Ends}
\DeclareMathOperator\QI{QI}
\DeclareMathOperator\interior{int}

\DeclareMathOperator\Stab{Stab}

\DeclareMathOperator\Hull{{\mathcal H}}

\DeclareMathOperator\Edges{{\mathcal E}}
\DeclareMathOperator\Verts{{\mathcal V}}
\DeclareMathOperator\PD{PD}
\DeclareMathOperator\Isom{Isom}
\DeclareMathOperator\kernel{Ker}

\newcommand\R{{\mathbf R}}

\newcommand\hyp{{\mathbf H}}

\newcommand\Z{{\mathbf Z}}
\newcommand\Euc{{\mathbf E}}
\renewcommand\H{{\mathcal H}}
\newcommand\F{{\mathcal F}}
\newcommand\E{{\mathcal E}}
\newcommand\C{{\mathcal C}}
\newcommand\fol\F
\renewcommand\P{{\mathbf P}}

\newcommand\solv{{\scshape solv}}

\newcommand\infinity{\infty}

\newcommand{\from}{\colon}

\newcommand\suchthat{\bigm|}
\newcommand\inverse{{-1}}
\newcommand\union{\cup}
\newcommand\Union{\bigcup}

\newcommand\intersect{\cap}

\newcommand\cross{\times}

\newcommand\Poincare{Poincar\'e}
\newcommand\cequiv{\approx_c}
\newcommand\ce\cequiv

\newcommand\cc{\subset_c}

\newcommand\ceq{=_c}
\newcommand\cintersect{\intersect_c}
\newcommand\<\langle
\renewcommand\>\rangle

\renewcommand\H{{\mathcal H}}

\newcommand\wt{\widetilde}

\begin{document}

\title{Quasi-actions on trees\\{\small Research announcement}}
\author{Lee Mosher, 
Michah Sageev, 
and Kevin Whyte}
\maketitle

We present a battery of techniques for investigating quasi-isometric
rigidity of graphs of groups \cite{MosherSageevWhyte:quasitree}. The
techniques work best when all edge and vertex groups are ``coarse
\Poincare\ duality'' groups in the sense of
\cite{KapovichKleiner:duality}, for example, fundamental groups of closed,
aspherical manifolds---such groups respond well to analysis using methods
of coarse algebraic topology introduced in \cite{FarbSchwartz} and further
developed in \cite{KapovichKleiner:duality}. Our techniques also require
the Bass-Serre tree to be bushy, meaning that it has infinitely many ends.
For example, when the dimensions of the vertex and edge groups are
homogeneous we obtain:

\begin{theorem}[Homogeneous Theorem] 
\label{ThmHomogeneous}
Fix an integer $n \ge 0$, and let
$\Gamma$ be a finite graph of coarse $\PD(n)$ groups with bushy 
Bass-Serre tree. Let $H$ be a finitely generated group quasi-isometric to
$\pi_1\Gamma$. Then $H$ is the fundamental group of a graph of groups with
bushy Bass-Serre tree, and with vertex and edge groups quasi-isometric to
those of $\Gamma$.
\end{theorem}

Specializing to $n=0$ gives a theorem of Stallings--Dunwoody
\cite{Stallings:ends}, \cite{Dunwoody:Accessible} which says, in modern
language, that if a finitely generated group
$H$ is quasi-isometric to a free group then $H$ is the fundamental group
of a finite graph of finite groups, implying furthermore that $H$ is
virtually free.

Theorem \ref{ThmHomogeneous} suggests the following problem. Given
$\Gamma$ as above, all edge--to--vertex group injections must have finite
index image, and so all vertex and edge groups lie in a single
quasi-isometry class. Given $\C$, a quasi-isometry class of coarse
$\PD(n)$ groups, let $\Gamma(\C)$ be the class of fundamental groups of
finite graphs of groups with vertex and edge groups in $\C$ and with bushy
Bass-Serre tree. Theorem \ref{ThmHomogeneous} says that 
$\Gamma(\C)$ is closed up to quasi-isometry.

\begin{problem} Given $\C$, describe the quasi-isometry classes within
$\Gamma(\C)$.
\end{problem}

Here's a rundown of the cases for which we know the solution to this
problem. In the case $n=1$, where $\C = \{\text{virtually
$\Z$}\}$: the amenable groups in $\Gamma(\C)$ are classified in
\cite{FarbMosher:BSOne}, and the nonamenable ones in \cite{Whyte:bs}. For
$\C = \{\text{virtually $\Z^n$}\}$, the amenable groups in $\Gamma(\C)$ are
classified in \cite{FarbMosher:abc}; the nonamenable case remains open. For
$\C =\{\text{quasi-isometric to $\hyp^2$}\} = \{\text{finite-by-(cocompact
fuchsian)}\}$, the subclass of $\Gamma(\C)$ consisting of word hyperbolic
surface-by-free groups is quasi-isometrically rigid and is classified in
\cite{FarbMosher:sbf}, but the broader classification in $\Gamma(\C)$ is
open. 

If $\C$ is the quasi-isometry class of cocompact lattices in an
irreducible, semisimple Lie group $L$ with finite center, then combining
Mostow Rigidity for $L$ with quasi-isometric rigidity (see
\cite{Farb:lattices} for a survey) it follows that for each $G \in \C$
there exists a homomorphism $G \to L$ with finite kernel and discrete,
cocompact image, and this homomorphism is unique up to post-composition
with an inner automorphism of $L$. Combining this with Theorem
\ref{ThmHomogeneous} it follows that any group in
$\Gamma(\C)$ is quasi-isometric to the cartesian product of any group in
$\C$ with any free group of rank $\ge 2$.

Our techniques also apply to graphs of coarse $\PD$ groups
without constant dimension, under various assumptions on how
edge spaces attach to the vertex spaces. Here is a sample application:

\begin{theorem}
\label{ThmAbelian} 
Let $\Gamma$ be a graph of finitely generated abelian
groups. Suppose that the following condition holds: 
\begin{itemize}
\item[(*)] for each vertex group $\Gamma_v$ of dimension $n$, all edge
groups incident to $\Gamma_v$ all have dimension $<n$, and if there exist
any incident edge groups of dimension $n-1$ then the incident edge groups
span $\Gamma_v$.
\end{itemize}
Then any finitely generated group $H$ quasi-isometric to $\Gamma$ is the
fundamental group of a graph of virtually abelian groups.
\end{theorem}

Sometimes we can combine our techniques with other results to get even
stronger quasi-isometric rigidity theorems. For example, by using the
Geodesic Pattern Rigidity Theorem of R.\ Schwartz
(\cite{Schwartz:Symmetric}, see Theorem \ref{ThmGeoPatterns}
below) we obtain:

\begin{theorem}
\label{ThmHVZE}
Let $\Gamma$ be a finite graph of groups whose vertex
groups are fundamental groups of closed hyperbolic manifolds of dimension
$\ge 3$, and whose edge groups are infinite cyclic. If $H$ is a finitely
generated group quasi-isometric to $\pi_1\Gamma$, then $H$ is the
fundamental group of a graph of groups $\Gamma'$, such that each vertex
group of $\Gamma'$ is weakly commensurable\footnote{\emph{Weak
commensurability} of groups $G_1,G_2$ means that there exists a group $Q$
and homomorphisms $Q \to G_1, Q \to G_2$ each with finite kernel and
finite index image. This is the smallest equivalence relation
incorporating passage both to a finite index subgroup and to a finite
kernel quotient.} either to a vertex group of $\Gamma$ or to $\Z$, and
such that each edge group of $\Gamma'$ is commensurable to $\Z$.
\end{theorem}

We will also find a way to strengthen Theorem \ref{ThmAbelian} by
developing Abelian Pattern Rigidity (Theorem \ref{ThmAbelPatterns}).

And there is more to the conclusions of Theorems \ref{ThmAbelian}
and \ref{ThmHVZE} than is stated: in each case the new graph of groups
not only has edge and vertex groups related to those of the old graph of
groups, but the edge-to-vertex injections are also related. This will be
made explicit in Theorem \ref{ThmInhom} and in the
applications of Theorem \ref{ThmInhom} to the proofs of
Theorems \ref{ThmAbelian} and \ref{ThmHVZE}.

\paragraph{Remark} Theorems \ref{ThmAbelian} and \ref{ThmHVZE} are
just samples. The main result on Inhomogeneous Rigidity, Theorem
\ref{ThmInhom}, has a multitude of applications. In this
research announcement we give only sketches of proofs, and we focus
particularly on coarse $\PD$ vertex and edge groups, ignoring wider
contexts for our results. Full statements in wider contexts,
and full proofs, will be found in \cite{MosherSageevWhyte:quasitree}.

\paragraph{Remark} As in the homogeneous case, the techniques for the
inhomogeneous case can be applied to quasi-isometric classification as
well as to quasi-isometric rigidity; although we have not mentioned here
any of these classification theorems, the preprint
\cite{MosherSageevWhyte:quasitree} will include some results. We should
also mention
\cite{PapasogluWhyte:ends}, where techniques similar to those of Theorem
\ref{ThmInhom} are used in proving that the quasi-isometric
classification of accessible groups completely reduces to the
classification of one-ended groups.

\section*{Graphs of groups and Bass-Serre trees: a review}

A graph of groups $\Gamma$ of finite type consists of the following data.
Start with a finite graph $\Gamma$ with vertex set $\Verts(\Gamma)$ and
edge set $\Edges(\Gamma)$. For each edge $e \in \Edges(\Gamma)$, the two
ends of $e$ form a set $\Ends(e)$, and each end $\eta \in \Ends(e)$ is
attached to some vertex $v(\eta) \in \Verts(\Gamma)$. Associated to each
vertex $v \in \Verts(\Gamma)$ there is a finitely generated vertex group
$\Gamma_v$, associated to each edge $e \in \Edges(\Gamma)$ there is an
edge group $\Gamma_e$, and associated to each end $\eta \in \Ends(e)$ there
is an edge-to-vertex injection $\gamma_\eta \from \Gamma_e \to
\Gamma_{v(\eta)}$. 

Associated to a graph of groups $\Gamma$ is a \emph{graph of spaces}
$B_\Gamma$, as follows. Choose path connected, pointed spaces
$B(v)$, $v\in \Verts(\Gamma)$, and $B(e)$, $e \in \Edges(\Gamma)$, whose
fundamental groups are identified with the associated vertex or edge groups
of $\Gamma$. Choose pointed \emph{attaching maps} $\xi_\eta\from B(e)\to
B(v)$ inducing the injections $\gamma_\eta$. For each edge $e$ let $\hat e
=
\interior(e) \union \Ends(e)$ denote the end compactification of $e$,
($\hat e$ is a compact arc, regardless of whether $e$ is a compact arc or
a loop). Construct a quotient space $B_\Gamma$ from the disjoint union of
the set
$$\{B(v),B(e) \cross \hat e \suchthat v\in \Verts(\Gamma), e \in
\Edges(\Gamma) \}
$$ 
by gluing $B(e) \cross \eta$ to $B(v(\eta))$ via the map $(x,\eta) \to
\xi_\eta(x)$, for each $e \in \Edges(\Gamma)$, $\eta \in\Ends(e)$. The
fundamental group $\pi_1(B_\Gamma)$ is well-defined up to isomorphism, and
it is called the fundamental group of $\Gamma$, denoted $\pi_1(\Gamma)$.

The map $B \to \Gamma$, taking $B(v)$ to $v$ and projecting $B(e) \cross
\hat e$ to $e$, induces a decomposition of $B$ into point
preimages. In the universal cover $X=\wt B$, taking connected lifts via
the universal covering map $X \to B$ of point preimages of $B\to\Gamma$
gives a decomposition of $X$ into path connected sets. This decomposition
of $X$ is $\pi_1(\Gamma)$-equivariant. The quotient space of this
decomposition of $X$ is a tree $T$ on which $\pi_1(\Gamma)$ acts,
the \emph{Bass-Serre tree} of $\Gamma$. This action is well-defined up to
equivariant tree isomorphisms, independent of the choices, and the
quotient graph of $T$ is canonically identified with $\Gamma$:
$$
\xymatrix{
X \ar[r] \ar[d] & T \ar[d] \\
B \ar[r] & \Gamma
}
$$
The map $X \to T$ is called the \emph{Bass-Serre tree of spaces}
associated to the graph of spaces $B\to\Gamma$. The inverse image of a
vertex $v \in T$ is a \emph{vertex space} $X(v)$ of $X$. The inverse image
of the midpoint of an edge $e$ of $T$ is called an \emph{edge space}
$X(e)$ of $X$. The topological space $X$ is constructed from the disjoint
union of the set $\{X(v), X(e) \cross e \suchthat v \in \Verts(T), e \in
\Edges(T)\}$, where for each vertex $v$ and edge $e$ incident to $v$ we
glue $X(e)\cross v$ to a subset of $X(v)$ via an attaching map $X(e) \to
X(v)$ which is a lift of an attaching map for the graph of spaces $B$. The
image of the attaching map $X(e) \to X(v)$ is called an \emph{incident
edge space inside~$X(v)$}. The set of incident edge spaces inside
$X(v)$ is called the \emph{edge space pattern inside $X(v)$}.

A simple trichotomy holds for Bass-Serre trees: $T$ is \emph{bounded};
or $T$ is \emph{line-like} meaning that it contains a line as a cobounded
subset; or $T$ is \emph{bushy} meaning that it has infinitely many ends. A
graph of groups $\Gamma$ (and its associated Bass-Serre tree) is said to
be \emph{reduced} if for each vertex $v$, the number of surjective
edge-to-vertex injections $\gamma_\eta \from\Gamma_e\to\Gamma_v$ is not
equal to $1$; it can be $0$ or $\ge 2$. When $T$ is reduced, the
trichotomy simplifies as follows. First, $T$ is bounded if and only if $T$
and
$\Gamma$ are each a point. Second, $T$ is line-like if and only if
$T$ is a line and $\Gamma$ is a \emph{mapping torus}, meaning either a
circle with isomorphic edge-to-vertex inclusions all around, or an arc
with isomorphic inclusions in the interior and index~2 inclusions at the
endpoints. Finally, $T$ is bushy if and only if $T$ has at least one
vertex of valence $\ge 3$; valence of a vertex in $T$ is easily computed
in terms of the image vertex in $\Gamma$, as the sum of the indices of the
edge groups inside the vertex group.

\section*{Coarse language}

Let $X$ be a metric space. Given $A \subset X$ and $R \ge 0$, denote
$N_R(A) = \{x \in X \suchthat \exists a \in A \quad\text{such that}\quad
d(a,x) \le R\}$. Given subsets $A, B \subset X$, let $A \cc B \,\,[R]$
denote $A \subset N_R(B)$. Let $A \cc B$ denote the existence of
$R \ge 0$ such that $A \cc B \,\,[R]$; this is called \emph{coarse
containment} of $A$ in $B$. Let $A \ceq B \,\,[R]$ denote $A \cc B \,\,[R]$
and $B\cc A \,\,[R]$. Let $A \ceq B$ denote the existence of $R$ such that
$A\ceq B \,\,[R]$; this is called \emph{coarse equivalence} of $A$ and $B$.

Given a metric space $X$ and subsets $A,B$, we say that a subset $C$ is a
\emph{coarse intersection} of $A$ and $B$ if for all sufficiently large
$R$ we have $N_R(A) \intersect N_R(B) \ceq C$. A coarse
intersection of $A$ and $B$ may not exist, but if one does exist then it
is well-defined up to coarse equivalence.

Given metric spaces $X,Y$, a map $f \from X \to Y$ is a \emph{uniform
embedding} if there exists proper, increasing functions $g,h \from
[0,\infinity)
\to [0,\infinity)$ such that
$$g(d_X(x,y)) \le d_Y(fx,fy) \le h(d_X(x,y))
$$
If $X,Y$ are geodesic metric spaces then the upper bound $h$ can always be
taken to be an affine function. When $h(d)=Kd+C$ and $g(d)=\frac{1}{K}d-C$
then we say that $f$ is a \emph{$K,C$ quasi-isometric embedding}. If this
is so then we say in addition that $f$ is a \emph{$K,C$ quasi-isometry
from $X$ to $Y$} if $f(X) \ceq Y \,\,[C]$. A \emph{$C'$-coarse inverse} of
$f$ is a $K,C'$ quasi-isometry $g \from Y \to X$ such that $x \ceq g(f(x))
\,\,[C']$ and $y \ceq f(g(y)) \,\,[C']$, for all $x \in X$, $y \in Y$. A
simple fact says that for all $K,C$ there exists $C'$ such that each $K,C$
quasi-isometry has a $C'$-coarse inverse.

Let $G$ be a group and $X$ a metric space. A \emph{$K,C$ quasi-action} of
$G$ on $X$ is a map $(g,x) \mapsto g \cdot x$ from $G \cross X$ to $X$,
such that: for each $g$ the map $x \mapsto g \cdot x$
is a $K,C$ quasi-isometry; and for each $x \in X, g,h \in G$ we have
$$g \cdot (h \cdot x) \ceq (gh) \cdot x \,\,[C]
$$
A quasi-action is \emph{cobounded} if there exists a constant $R$ such that
for each $x \in X$ we have $G \cdot x \ceq X \,\,[R]$. A quasi-action
is \emph{proper} if for each $R$ there exists $M$ such that for all $x,y
\in X$, the cardinality of the set $\{g \in G \suchthat \bigl(g \cdot
N(x,R)\bigr) \intersect N(y,R) \ne \emptyset\}$ is at most $M$.

A fundamental principle of geometric group theory says that if $G$ is a
finitely generated group equipped with the word metric, and if $X$ is a
proper geodesic metric space on which $G$ acts properly discontinuously
and cocompactly by isometries, then $G$ is quasi-isometric to $X$.

A partial converse to this result is the \emph{quasi-action principle}
which says that if $G$ is a finitely generated group with the word metric
and $X$ is a metric space quasi-isometric to $G$ then there is a
cobounded, proper quasi-action of $G$ on $X$; the constants for this
action depend only on the quasi-isometry constants between $G$ and
$X$. 

The quasi-action principle motivates the following question, which is a
common point of departure for many quasi-isometric rigidity problems:
\begin{itemize}
\item Given a proper, cobounded quasi-action of a group $G$ on a metric
space $X$, when can we get some action of $G$ on $X$? 
\end{itemize}
Partial information about this question can sometimes be obtained from the
following result (see e.g.\ \cite{KapovichLeeb:haken}):

\begin{proposition}[Coboundedness Principle]
\label{PropCobounded} 
Suppose that a finitely generated $G$ quasi-acts
properly and coboundedly on a metric space $X$. Let $\H$ be a collection
of subsets of $X$ which satisfies the following properties: the elements
of $\H$ are pairwise coarsely inequivalent in $X$; there exists $A \ge 0$
such that for each
$g \in G$, $H \in \H$ there is an $H' \in\H$ such that $g \cdot H \ceq H'
\,\,[A]$; and every metric ball in $X$ intersects at most finitely many
elements of $\H$. Then for each $H \in\H$, the \emph{stabilizer subgroup}
$\Stab_G(H) = \{g \in G \suchthat g\cdot H \ceq H\}$ quasi-acts properly
and coboundedly on $H$.
\end{proposition}

We will need some coarse algebraic topology. This subject originated in
\cite{FarbSchwartz}, by applying Alexander Duality with naturality to
fundamental groups of closed aspherical manifolds. These ideas were
generalized in \cite{KapovichKleiner:duality} to work for a general
class of spaces, the \emph{coarse $\PD(n)$ spaces $X$}. These are
``uniformly acyclic'' simplicial complexes which satisfy a coarse, uniform
version of the \Poincare\ duality property of $\R^n$: there exist chain
maps $C_*(X)
\xrightarrow{D} C^{n-*}_c(X)$ and $C^{n-*}_c(X)
\xrightarrow{\overline D} C_*(X)$ with ``uniform distortion'' such that $D
\circ \overline D$ and $\overline D \circ D$ are ``uniformly chain
homotopic'' to the identity; see \cite{KapovichKleiner:duality} for the
full development. A coarse $\PD(n)$ group $G$ is one which acts properly
and coboundedly on some coarse $\PD(n)$ space. 

A subset $H$ of a metric space $X$ is \emph{deep} if for each $r$ there
exists $x \in X$ such that the ball in $X$ around $x$ of radius $r$ is
a subset of $H$.

\begin{proposition}[Coarse Jordan Separation] 
\label{PropCoarseJordan}
If $X$ is a coarse $\PD(n)$
space and $S\subset X$ is a uniformly embedded coarse $\PD(n-1)$ space,
then for sufficiently large $A$ there are exactly two deep components of
$X-N_A(S)$. The coarse intersection of these two components is $S$.
\end{proposition}

For fundamental groups of closed, aspherical $n$-manifolds, this result
with a weaker conclusion of ``at least two deep components'' comes from
\cite{FarbSchwartz}. The improvement to ``exactly two deep components'',
and the generalization to coarse $\PD(n)$ spaces, comes from
\cite{KapovichKleiner:duality}.

We will state other coarse algebraic topology results as we need them
below.

\section*{Quasi-actions on trees}

Suppose that $\Gamma$ is a finite type graph of groups, $B \to \Gamma$ is
an associated graph of spaces, and $X \to T$ is the associated Bass-Serre
tree of spaces. We make the additional assumption that each edge and
vertex space of $B$ is compact; for instance, if all edge and vertex
groups are finitely presented then we can take the edge and vertex spaces
to be presentation complexes. We impose a geodesic metric on $B$, which
lifts to a geodesic metric on $X$. It follows that $\pi_1\Gamma$ is
quasi-isometric to $X$. If $G$ is any finitely generated group
quasi-isometric to $\pi_1\Gamma$, it follows that $G$ is also
quasi-isometric to $X$, and therefore $G$ has a cobounded, proper
quasi-action on $X$. In this situation, the motivating question is:
\begin{itemize}
\item Does the quasi-action of $G$ on $X$ coarsely respect the vertex and
edge spaces of $X$?
\end{itemize}
While this question is somewhat vague, there are several distinct ways to
make it more precise. Here is the most rigid possible behavior:
\begin{itemize}
\item \textbf{Tree rigidity\/} There is an action of $G$ on $T$ such that
for each $g\in G$, $v\in\Verts(T)$, and $e \in \Edges(T)$ we have $g \cdot
X(v)
\ceq X(g\cdot v)$ and $g \cdot X(e) \ceq X(g \cdot e)$ (uniformly,
i.e.\ with uniform coarseness constant independent of $g,v$).
\end{itemize}
Tree rigidity immediately implies, for example, that $G$ is the
fundamental group of a graph of groups whose vertex and edge groups are
all quasi-isometric to vertex and edge groups of $\Gamma$. 

In general, tree rigidity is a lot to expect. Ignoring edges for the
moment, here is a sequence of vertex rigidity properties, from weaker to
stronger:
\begin{itemize}
\item \textbf{Weak vertex rigidity\/} For each $g \in G$, $v \in \Verts(T)$
there exists $w \in \Verts(T)$ such that $g \cdot X(v) \cc X(w)$
[uniformly].
\item \textbf{Vertex rigidity\/} For each $g \in G$, $v \in \Verts(T)$
there exists $w \in \Verts(T)$ such that $g \cdot X(v) \ceq X(w)$
[uniformly].
\item \textbf{Vertex rigidity with uniqueness\/} For each $g \in G$, $v
\in \Verts(T)$ there exists a unique $w \in \Verts(T)$ such that $g \cdot
X(v) \ceq X(w)$ [uniformly].
\end{itemize}
In general these three properties are not equivalent. In certain
situations they become equivalent, e.g.\ when no vertex space is
coarsely contained in another, which occurs if and only if every
edge-to-vertex group injection has infinite index (a property which fails
completely in the context of Theorem \ref{ThmHomogeneous}).

The utility of vertex rigidity with uniqueness, for instance, is that it
allows us to define an action of $G$ on $\Verts(T)$, where $g \cdot v = w$
if and only if $g \cdot X(v) \ceq X(w)$. Even then we don't get any action
of $G$ on the edges, without knowing something more about edge spaces,
such as:
\begin{itemize}
\item \textbf{Strong edge rigidity} For all edges $e \ne e' \in
\Edges(T)$, $X(e_ \not\ceq X(e')$.
\end{itemize}
Vertex rigidity with uniqueness, coupled with strong edge rigidity,
implies tree rigidity, for if $e$ is an edge with endpoints $v_1,v_2$, and
if $e'$ is the edge with endpoints $g \cdot v_1, g \cdot v_2$, then we have
$$g \cdot X(e) \ceq g \cdot (X(v_1) \cintersect X(v_2)) \ceq g \cdot
X(v_1) \cintersect g \cdot X(v_2) \ceq X(g \cdot v_1) \cintersect X(g
\cdot v_2) \ceq X(e')
$$
and so by setting $g \cdot e = e'$ we obtain a well-defined action of $G$
on $T$ satisfying tree rigidity. 

\paragraph{Example} By the work of Kapovich--Leeb
\cite{KapovichLeeb:haken}, the torus decomposition of a non-\solv\ Haken
\nb{3}manifold satisfies tree rigidity. The heart of their work is a proof
of weak vertex rigidity, using asymptotic cones. Tree rigidity follows
easily from that, because the torus decomposition has infinite index
edge-to-vertex group injections, and obviously satisfies strong edge
rigidity.

\section*{Homogeneous graphs of groups: proof of Theorem
\protect\ref{ThmHomogeneous}}

A graph of groups is \emph{geometrically homogeneous} if it it satisfies
any of the following equivalent statements: each edge-to-vertex injection
is a quasi-isometry; each edge-to-vertex injection has finite index image;
the Bass-Serre tree has bounded valence. It follows that the edge and
vertex groups are all quasi-isometric. 

Here's the first step in the proof of Theorem \ref{ThmHomogeneous}:

\begin{theorem} 
\label{ThmFM}
Suppose $\Gamma,\Gamma'$ are geometrically homogeneous
graphs of coarse $\PD$ groups with bushy Bass-Serre trees $T,T'$, and
trees of spaces $X,X'$. If $h \from X \to X'$ is a quasi-isometry then $h$
respects vertex spaces. More precisely, for each $K,C$ there exists
$A$ such that if $h \from X \to X'$ is a $K,C$ quasi-isometry then
for each $v\in \Verts(T)$ there exists $v' \in \Verts(T')$ such that
$h(X(v)) \ceq X'(v') \,\,[A]$.
\end{theorem}

A proof with some extra topological assumptions is found in
\cite{FarbMosher:abc}, and that proof generalizes almost word-for-word to
coarse $\PD(n)$ vertex and edge groups. Applying Theorem
\ref{ThmFM} to each element of a quasi-action we obtain:

\begin{corollary} 
\label{CorTreeQuasi}
With the same notation as in Theorem
\ref{ThmFM}, if $H$ is a finitely generated group
quasi-isometric to $\pi_1\Gamma$ then the quasi-action of $H$ on $X$
satisfies vertex rigidity. It follows that there is a cobounded
quasi-action of $H$ on $\Verts(T)$, with the property that
$h \cdot X(v) \ceq X(h \cdot v)$ [uniformly].
\end{corollary}

The following result is the technical heart of Theorem
\ref{ThmHomogeneous}:

\begin{theorem}
\label{ThmBddValence}
Let $T$ be a bushy tree of uniformly bounded valence. Let $H$ be a group
quasi-acting coboundedly on $T$. Then the quasi-action of $H$ on $T$ is
quasiconjugate to a cobounded action of $H$ on a tree $T'$. That is, there
is a quasi-isometry $f\from T \to T'$ such
that $h \cdot (f(v)) \ceq f(h \cdot v)$ [uniformly].
\end{theorem}  

This combines with Corollary \ref{CorTreeQuasi} to prove
Theorem \ref{ThmHomogeneous}, for one can show first that the tree $T'$ has
bounded valence, and so there is a graph of groups
$\Gamma'$ with fundamental group $H$ and Bass-Serre tree $T'$.
Furthermore, the quasi-isometry $H \to \pi_1(\Gamma)$ takes the vertex and
edge stabilizers of the action of $H$ on $T'$ to subsets of $\pi_1\Gamma$
coarsely equivalent to the vertex and edge stabilizers of the action of
$\pi_1\Gamma$ on $T$.

\begin{proof}[Proof of Theorem \protect\ref{ThmBddValence}]
The first step of the proof is to construct a vertex set and an
action of $H$. We do not even have a true action of $H$ on
$\Verts(T)$, however. What we do have is an action of $H$ on the space of
ends $\Ends(T)$. This can be promoted into an action of $H$ on a new
vertex space as follows. Bushiness and bounded valence of $T$ tells us that
the space of ends
$\Ends(T)$ is a Cantor set. Note that each edge of $T$ determines a
partition of $\Ends(T)$ into a pairwise disjoint set of clopens
(closed-open subsets). Define a \emph{quasi-edge} of $T$ to be any
partition of $\Ends(T)$ into a pair of disjoint clopens. Clearly $H$ acts
on the set of quasi-edges of $T$. If $O \subset \Ends(T)$ is a clopen then
the convex hull $\Hull(O)$ is a subtree of $T$. If $\{O,O'\}$ is a
quasi-edge then the intersection of the \nb{1}neighborhoods of the convex
hulls  $N_1(\Hull(O)) \intersect N_1(\Hull(O'))$ is bounded; the diameter
of this intersection is defined to be the \emph{quasi-edge constant} of
$\{O,O'\}$. Note that a \nb{1}quasi-edge is the same thing as (the
partition determined by) a true edge of $T$. The action of a
$(K,C)$ quasi-isometry of $T$ on an $A$-quasi-edge produces an
$A'$-quasi-edge, where $A'$ depends only on $K,C,A$. It follows that
the $H$-orbit of a single quasi-edge has uniform quasi-edge constant. 

We therefore \emph{define} a vertex set $Y^0$ to be the $H$-orbit of some
arbitrarily chosen quasi-edge of $T$. It is not hard, using uniformity of
the quasi-edge constant, to attach edges to $Y^0$ in a locally finite,
$H$-equivariant manner, thereby producing a locally finite graph $Y^1$ on
which $H$ acts coboundedly, and a coarse conjugation $T\to Y^1$.

Unfortunately, $Y^1$ may not be a tree, i.e.\ in particular it may not be
simply connected. The latter problem can be corrected without too much
difficulty by attaching \nb{2}cells to $Y^1$ in a locally finite,
$H$-equivariant manner, producing a locally finite, simply connected
\nb{2}complex $Y^2$ on which $H$ acts coboundedly; the coarse conjugation
$T\to Y^1$ extends to a coarse conjugation $T \to Y^2$.

The final step is to get from the \nb{2}complex $Y^2$ back to a tree. This
is accomplished by using tracks in the sense of Dunwoody
\cite{Dunwoody:Accessible}: a track $\tau$ in $Y^2$ is a locally
separating, connected \nb{1}complex in general position with respect to
the skeleta of $Y^2$. Since $Y^2$ is simply connected, each track $\tau$
separates $Y^2$ into two components, and $\tau$ is \emph{essential} if
each component of $Y-\tau$ is unbounded. Using minimal surface ideas as in
\cite{Dunwoody:Accessible}, we construct an $H$-equivariant family of
pairwise disjoint, essential tracks $\{\tau_i\}$ in $Y^2$. By a Haken
finiteness argument as in \cite{Dunwoody:Accessible} one shows that a
maximal such family
$\{\tau_i\}$ exists, and that all components of $Y - \Union_i \tau_i$ are
bounded. The final tree $T'$ has vertex set in 1--1 correspondence with
the components of $Y-\Union_i \tau_i$, and edge set in 1--1 correspondence
with the set $\{\tau_i\}$.
\end{proof}

Theorem \ref{ThmBddValence} has an interesting consequence
concerning lattices in locally compact topological groups. Some
quasi-isometry classes $\C$ have the nice property that there is a locally
compact group $H$ such that each group $G \in \C$ has a discrete,
cocompact quotient in $H$ with finite kernel. This is true, for example,
when $\C$ is the quasi-isometry class of cocompact lattices in a
semisimple Lie group. However, this doesn't work for $\C =
\{\text{virtually free of rank $\ge 2$}\}$:

\begin{corollary}
There does not exist a locally compact group $H$ such that every
group which is virtually free of rank $\ge 2$ has a discrete, cocompact
quotient in $H$ with finite kernel.
\end{corollary}

\begin{proof} Suppose $H$ exists as stated. Choose a free, discrete,
cocompact subgroup $G \subset H$. We may give $H$ a left-invariant metric
so that the inclusion of $G$ into $H$ is a quasi-isometry. Since
$H$ is quasi-isometric to the free group $G$, it follows that $H$ is
quasi-isometric to a finite valence tree $T$. The left action of $H$ on
itself is quasi-conjugate to a cobounded quasi-action of $H$ on $T$.
Applying Theorem
\ref{ThmBddValence} it follows that $H$ has a cobounded
quasi-action on a bushy tree $T'$ of bounded valence. It follows that
\emph{every} virtually free group $G$ has a cobounded action on $T'$. Pick
a prime $p$ larger than the maximal valence of a vertex of $T'$. The group
$G = \Z/p * \Z/p$ is virtually free, and therefore has discrete cocompact
image in $H$ with finite kernel, and so $G$ acts coboundedly on $T'$. But
each of the free factors $\Z/p$ acts trivially on $T'$, making $G$ act
trivially, contradicting coboundedness.
\end{proof}

\section*{Inhomogeneous graphs of groups}

Theorem \ref{ThmAbelian}, and part of the conclusion of Theorem
\ref{ThmHVZE}, will both follow from more general results about
geometrically inhomogeneous graphs of groups. Although as in the
homogeneous case we have results in wider contexts, for present purposes
we restrict to coarse $\PD(n)$ vertex and edge groups. 

Suppose that $\Gamma$ is a graph of coarse $\PD(n)$ groups of various
dimensions. An \emph{$n$-raft} is a connected subgraph of $\Gamma$ (or of
$T$) of constant dimension $n$, such that each edge incident to the raft
but not contained in the raft has dimension $<n$. Rafts in $T$ are
connected lifts of rafts in $\Gamma$. Each raft in $T$ is the Bass-Serre
tree for the corresponding quotient raft in $\Gamma$. 

Assuming $\Gamma$ is reduced, each line-like raft in $T$ is actually a
line, and the quotient is a mapping torus raft in $\Gamma$. The problem
with mapping tori is that they generally fail to satisfy even weak vertex
rigidity. For example, a closed hyperbolic 3-manifold fibering over $S^1$
has a fundamental group whose quasi-isometry group is all of $\QI(\hyp^3)
= QC(S^2)$, and weak vertex rigidity fails miserably. A lattice in
\nb{3}dimensional \solv\ geometry is a mapping torus of $\Z^2$, and while
vertex rigidity is conjectured to hold \cite{FarbMosher:solvable}, the
conjecture remains open. For these and other reasons, in all of our
further theorems we must avoid line-like rafts in Bass-Serre trees.

Consider an $n$-dimensional point raft $v$ of the Bass-Serre tree $T$ with
associated vertex space $X(v) \subset X$. Let $\E$ be the edge space
pattern inside $X(v)$. For each $m < n$ define a subset $\E_m \subset
\E$ of $m$-dimensional edge spaces inside
$X(v)$, and let $\E_{[m,m']} =\E_m \union \E_{m+1}\union\cdots\union
\E_{m'}$. The Coarse Jordan Separation Theorem implies that each $E\in
\E_{n-1}$ coarsely separates $X(v)$ into two deep pieces, whose coarse
intersection is $E$; in this situation we say that a subset of $X(v)$
\emph{crosses} $E$ if it intersects each of the two deep pieces arbitrarily
far from $E$. Define the \emph{crossing graph} of
$X(v)$ to be the graph whose vertex set is $\E_{n-1}$, with an edge between
$E,E' \in \E_{n-1}$ if $E,E'$ cross each other or if there exists some
element of $\E_{[1,n-2]}$ which crosses both $E$ and $E'$.

\begin{theorem}
\label{ThmInhom}
Let $\Gamma$ be a finite, reduced graph of coarse $\PD$ groups satisfying
the following \emph{``raft hypotheses''}: the Bass-Serre tree $T$ has no
line rafts; and for each point raft $v$ of $T$, the crossing graph of
$X(v)$ is connected or empty. If $H$ is a finitely generated group
quasi-isometric to $\pi_1\Gamma$, then there is a finite type, reduced
graph of groups $\Gamma'$ with $H\approx\pi_1\Gamma'$ and with Bass-Serre
tree of spaces $X'\to T'$, and there is a quasi-isometry $f\from X' \to X$
coarsely conjugating the $H$ action on $X'$ to the $H$-quasi-action on $X$,
such that $f$ coarsely respects rafts and their vertex spaces, and $f$
coarsely respects all edge spaces. That is: for each raft $t'$ of $X'$
there exists a raft $t$ of $X$ such that $f(X'(t')) \ceq X(t)$
[uniformly], for each vertex $v' \in t'$ there exists a vertex $v \in t$
such that $f(X'(v')) \ceq X(v)$ [uniformly]; and for each edge $e'$ of $T'$
there exists an edge $e$ of $T$ such that $f(X'(e')) \ceq X(e)$
[uniformly].
\end{theorem}

Although the conclusion makes no mention of vertices $v' \in T'$ not on
any raft of $T'$, such information may be derived as follows: since
$v'$ is not on any raft, there exists an edge $e' \subset T'$ incident to
$v'$ such that $X'(v') \ceq X'(e')$, and so $f(X'(v')) \ceq f(X'(e'))$ is
coarsely equivalent to some edge space in $X$. However, counterexamples
show that $f(X'(v'))$ might not be coarsely equivalent to any vertex space
in $X$.

Before proving the theorem, we apply it to Theorem \ref{ThmAbelian}
and part of Theorem \ref{ThmHVZE}. 

The hypothesis (*) in Theorem \ref{ThmAbelian} immediately implies the
raft hypotheses in Theorem \ref{ThmInhom}, and the conclusion
of Theorem \ref{ThmInhom} immediately implies the conclusion
of Theorem \ref{ThmAbelian}. Note however that Theorem
\ref{ThmInhom} gives a much stronger conclusion, namely that
the new Bass-Serre tree of spaces $X'\to T'$ maps to original Bass-Serre
tree of spaces $X \to T$, preserving the coarse inclusions among vertex
and edge spaces. We'll use this additional information later, combining it
with our Abelian Pattern Rigidity Theorem \ref{ThmAbelPatterns} to
obtain a stronger rigidity result.

For Theorem \ref{ThmHVZE}, the fact that vertex groups have dimension
$\ge 3$ and edge groups have dimension $1$ implies that all rafts are
point rafts with empty crossing graph, and so Theorem
\ref{ThmInhom} applies. We conclude that all edge groups of
$\Gamma'$ are quasi-isometric to $\Z$, and therefore commensurable to
$\Z$. We also conclude that each vertex group $\Gamma'_{v'}$ satisfies one
of two possibilities: $\Gamma'_{v'}$ is quasi-isometric and so
commensurable to $\Z$; or $\Gamma'_{v'}$ quasi-isometric to some vertex
group $\Gamma_v$, and so to $\hyp^n$, and so $\Gamma'_{v'}$ is weakly
commensurable to \emph{some} closed $\hyp^n$ orbifold group.  But in the
latter case we get more information: the ambient quasi-isometry
$X'\to X$ takes $X'(v')$ to $X(v)$, coarsely mapping the edge
space pattern inside $X'(v')$ to the edge space pattern inside $X(v)$. 
We'll make this more precise later, and combine it with Schwartz' Geodesic
Pattern Rigidity Theorem \ref{ThmGeoPatterns}, to get the stronger
conclusion of Theorem \ref{ThmHVZE}, namely that $\Gamma'_{v'}$ is
weakly commensurable to $\Gamma_v$ itself.

\section*{Sketch of proof of Theorem \protect\ref{ThmInhom}}

Let $N$ be the maximal dimension of a vertex in $T$, and define a
filtration $T_N \subset\cdots\subset T_{i+1} \subset T_i \subset \cdots
\subset T_0=T$ where $T_i$ is the union of all vertices and edges of
dimension $\ge i$. Note that $T_N$ is a disjoint union of $N$-rafts. There
may be lower dimensional rafts as well: any component of $T_i$ which does
not contain a component of $T_{i+1}$ is an $i$-raft. Let $\Gamma_i$
be the image of $T_i$ in $\Gamma$, and let $X_i$ be the inverse image of
$T_i$ in $X$. Thus, each component of $X_i$ maps to a component of
$T_i$, giving a Bass-Serre tree of spaces for the corresponding
component of $\Gamma_i$.

Let $H$ be a group quasi-isometric to $\pi_1\Gamma$, and so $H$ quasi-acts
properly and coboundedly on $X$. The method of the proof is to work
inductively down from the top dimension, replacing the quasi-action of $H$
on the tree of spaces $X_i \to T_i$ by a true action, starting with $i=N$. 
The basis step of the induction depends on the following:

\begin{proposition}[Vertex rigidity at the top dimension] 
\label{PropTopRigidity}
The quasi-action
of $H$ on $X$ coarsely respects $N$-rafts and their vertex spaces, that
is: for each $N$-raft $t \subset T$ and each $h\in H$ there exists an
$N$-raft $t' \subset T$ such that $h \cdot X(t) \ceq X(t')$ [uniformly],
and for each vertex $v \in t$ there exists a vertex $v'\in t'$ such that $h
\cdot X(v) \ceq X(v')$ [uniformly]. 
\end{proposition}

Once this is shown, applying Proposition \ref{PropCobounded} to the
collection of $N$-rafts we conclude that the stabilizer of each $N$-raft
quasi-acts coboundedly and properly, and then applying Theorem
\ref{ThmHomogeneous} we may quasi-conjugate the action of $\Stab(T)$
on each $N$-raft $T$ to a true action of $\Stab(T)$ on a new tree $T'$.
Thus we establish the basis step for the inductive proof of Theorem
\ref{ThmInhom}.

\begin{proof}[Proof of Proposition \protect\ref{PropTopRigidity}] We shall
prove: 
\begin{itemize}
\item[(*)] There exists $A$ such that for each $v \in T_N$ and each $h \in
H$ there is an $N$-raft $t$ with $h \cdot X(v) \cc X'(t)
\,\,[A]$. 
\end{itemize}
To see why this suffices, consider a raft $t$ of $T_N$ and two vertices
$v_1,v_2 \in \Verts(t)$, and so $X(v_1) \ceq X(v_2)$. Applying $(*)$ we
get $h\cdot X(v_1)\cc X'(t_1)$, $h \cdot X(v_2) \cc X'(t_2)$ for
$N$-rafts $t_1,t_2$ of $T$, and we want to verify that $t_1=t_2$. Since
$X(v_1) \ceq X(v_2)$ it follows that the coarse intersection of $X'(t_1)$
and $X'(t_2)$ coarsely contains a coarse $\PD(N)$ space, namely $h \cdot
X(v_1)\ceq h \cdot X(v_2)$. If
$t_1\ne t_2$ then the coarse intersection of $X'(t_1)$ and $X'(t_2)$ is
coarsely equivalent to some edge space of lower dimension, but a coarse
$\PD(n)$ space with $n<N$ cannot contain a uniformly embedded copy of a
coarse $\PD(N)$ space. It follows that $t_1=t_2$ and $h \cdot X(t) \cc
X(t_1)$. Similarly $h^\inv \cdot X(t_1) \cc X(t')$ for some $N$-raft $t'$,
and so $X(t') \cc X(t)$; but this is only possible if $t=t'$.
This shows that $H$ coarsely respects $N$-rafts. 

If (*) is not true then, taking counterexamples for larger and larger
values of $A$, using coboundedness of the isometry group of $X\to T$, and
passing to a limit, we obtain a quasi-isometry $h \from X \to X$, a
vertex $v\in\Verts(T_N)$ and an edge $e \in\Edges(T)-\Edges(T_N)$, such
that
$h( X(v))$ intersects both components of $X-X(e)$ arbitrarily
deeply. It follows that for sufficiently large $A$, the subset
$$S = h^\inv(N_A(X(e))) \intersect X(v)
$$
coarsely separates $X(v)$ into at least two deep components. But the set
$S$, with metric restricted from $X(v)$, is uniformly equivalent to a
subset of $X(e)$ with restricted metric. 

When $X(e)$ has dimension $\le N-2$ we obtain a contradiction
using arguments of coarse algebraic topology: a coarse $\PD(N)$ space
cannot be coarsely separated by a subset which is uniformly equivalent to
a subset of a coarse $\PD$ space of dimension $ \le N-2$. 

If $X(e)$ is of dimension $=N-1$ then in fact $\pi(h(S)) \ceq X(e)$:
otherwise, a subset of the coarse $\PD(N-1)$ space $X(e)$ which is not
coarsely equivalent to all of $X(e)$ would embed uniformly in coarse
$\PD(N)$ space $X(v)$, coarsely separating $X(v)$, and that is impossible.
We therefore have $h(S) \ceq X(e)$ in $X$. This shows that $S$ is a coarse
$\PD(N-1)$ space uniformly embedded in the coarse $\PD(N)$ space $X(v)$,
and so $X(v) - S$ has exactly two deep components each coarsely containing
$S$, each contained in a separate deep component of $X-S$. 

Now the argument breaks into cases.

\textbf{Case 1: $v$ is contained in a bushy raft $t$.} In this case
the coarse $\PD(N-1)$ space $S$ coarsely separates $X(t)$, a tree of
$\PD(N)$ spaces, and that is clearly impossible.

\textbf{Case 2: $v$ is a point raft.} By hypothesis, the
crossing graph of $X(v)$ is either connected or empty.

\textbf{Case 2a: The crossing graph of $X(v)$ is connected.} Using
connectedness of the crossing graph together with some coarse separation
arguments, one shows that one of the deep components of $X(v) - S$ coarsely
contains the union of all codimension-1 edge spaces inside $X(v)$. But this
is absurd, because coboundedness of the $\pi_1\Gamma$ stabilizer subgroup
of $X(v)$ shows that the union of incident codimension-1 edge spaces
intersects every deep subset of $X(v)$.

\textbf{Case 2b: The crossing graph of $X(v)$ is empty.} Each
edge incident to $v$ therefore has dimension $\le N-2$, and it follows
that the inclusion of $S$ in $X$ has the following ``coarse Jordan
separation property'': for all sufficiently large $A\ge 0$ there are
\emph{exactly} two deep components of $X-N_A(S)$ which coarsely contain
$S$. The inclusion of $X(e) \ceq h(S)$ into $X$ therefore has the same
property: for all sufficiently large $A$ there are exactly two deep
components of $X-N_A(X(e))$ which coarsely contain $X(e)$. 

Let $T^*_e$ be the subtree of $T$ spanned by all edges in $T$ whose edge
space is coarsely equivalent to $X(e)$; we think of $T^*_e$ as the ``edge
raft'' containing $e$, although a priori $T^*_e$ can have $N$-dimensional
vertices and edges. But by using the coarse Jordan separation property for
$X(e)$ one shows that $T^*_e$ contains at most one $N$-dimensional vertex of
$T$. Moreover, $T^*_e$ cannot contain exactly one $N$-dimensional vertex,
for then one would be able to find an $N-1$ dimensional valence~1 vertex of
$T^*_e$ which would violate irreducibility of the Bass-Serre tree $T$. It
follows that $T^*_e$ in fact consists entirely of $N-1$ dimensional
vertices and edges, and so is an $N-1$ raft. $T^*_e$ cannot be a bounded
raft, for it has at least one edge, namely $e$, and again that would
violate irreducibility. $T^*_e$ cannot be a line raft, by hypothesis.
Finally, it cannot be a bushy raft, because that would violate the coarse
Jordan separation property for $X(e)$: for larger and larger $A$, the
number of deep components of $X-N_A(X(e))$ coarsely containing $X(e)$ would
approach infinity. 
\end{proof}

We now continue the proof of Theorem \ref{ThmInhom} by
induction down the dimension. Suppose by induction that we have altered
$X$ and $T$ down to dimension $n+1$, producing a filtered tree of spaces
$$
\xymatrix{
X'_N \ar[d] \ar[r]^\subset & \cdots \ar[r]^\subset & X'_{n+1} \ar[d]
\ar[r]^\subset & X_n \ar[d]
\ar[r]^\subset & \cdots \ar[r]^\subset & X_0 \ar[d] \\
T'_N \ar[r]_\subset & \cdots \ar[r]_\subset & T'_{n+1} \ar[r]_\subset & T_n
\ar[r]_\subset & \cdots \ar[r]_\subset & T_0 }
$$
so that $H$ quasi-acts properly and coboundedly on $X_0$, restricting
to a true action on $X'_k \to T'_k$, $N \ge k \ge n+1$, and we have an
$H$-quasiconjugation back to the original tree of spaces, restricting to
the identity on $X_0 - X'_{n+1}$, and satisfying the conclusions of
Theorem \ref{ThmInhom} on $X'_{n+1}$.

The arguments of the basis step can be applied to any component of $T_n$
which is an $n$-raft.

Consider now an edge $e$ of $T_n - T'_{n+1}$ which does not lie on an
$n$-raft, and an element $h \in H$; we study the image $h \cdot X(e)$.
Using irreducibility of the original tree of spaces it follows that there
are two vertices $v,w \in T'_{n+1}$ such that $X(e)$ is a coarse
intersection of $X'(v)$ and $X'(w)$, with constants independent of $e$;
moreover, $v,w$ can be chosen to have a distance bounded above in
$T'_{n+1}$ independent of $e$. Now we apply the inductive hypothesis to $h
\cdot X(v), h \cdot X(w)$, splitting into subcases depending on whether
$v,w$ lie on rafts. 

\medskip

\textbf{Case 1:} Suppose $v,w$ do lie on rafts. By induction, there exist
vertices $v',w'$ such $h \cdot X(v) \ceq X(v')$, $h\cdot X(w) \ceq X(w')$.
It follows that $h \cdot X(e)$ is a coarse intersection of $X(v')$ and
$X(w')$. Let $e_1,\ldots,e_K$ be the simple edge path in $T_0$ connecting
$v'$ to $w'$; $K$ is bounded independent of $v,w$, depending only on the
quasi-isometry constants of $h$. The coarse intersection of
$X(v')$ and $X(w')$ equals the coarse intersection of
$X(e_1),\ldots,X(e_K)$. This implies that $h\cdot X(e)$ is a coarse
intersection of $X(e_1),\ldots,X(e_K)$, and in particular $h
\cdot X(e)$ is coarsely contained in each of $X(e_1),\ldots,X(e_K)$. 

The first consequence of this is that the dimensions of $e_1,\ldots,e_K$
are greater than or equal to the dimension of $e$, because a coarse
$\PD(n)$ space cannot uniformly embed in a coarse $\PD$ space of lower
dimension. 

The second consequence is that if $e_k$ is $n$-dimensional then $h \cdot
X(e) \ceq X(e_k)$, because of ``packing'': a uniform embedding of a coarse
$\PD(n)$ space in a coarse $\PD(n)$ space must have image coarsely
equivalent to the whole space.

To complete the proof in Case 1 it therefore remains to check that the
edges $e_1,\ldots,e_K$ cannot all have dimension $\ge N+1$. For this we
need the fact that $X(e)$ coarsely separates $X(v)$ and $X(w)$ in $X$, and
it follows that $h \cdot X(e)$ coarsely separatex $X(v')$ and $X(w')$ in
$X$. But if $e_1,\ldots,e_K$ all have dimension $\ge N+1$ then the coarse
$\PD(n)$ space $h \cdot X(e)$ cannot coarsely separate $X(v')$ and $X(w')$.

\medskip

\textbf{Remaining cases:} If, say, $v$ lies on a raft and $w$ does not,
then $w$ is incident to an edge $e_1$ of the same dimension. Applying
induction, $h \cdot X(v) \ceq X(v')$ for some raft vertex $v'$, and $h
\cdot X(w) \ceq h \cdot X(e_1) \ceq X(e')$ for some edge $e'$. Now connect
$v'$ and $e'$ by a simple edge path and repeat the arguments of Case 1.
The other cases are similar.

\medskip

We have shown that the quasi-action of $H$ coarsely respects
$n$-dimensional edge spaces. It remains to attach additional edges to the
$H$-forest $T'_{n+1}$ to make an $H$-forest $T'_n$, so that a newly
attached edge $e'$ between vertices $v', w' \in T'_{n+1}$ has an coarse
$\PD(n)$ edge space $X'(e') \subset X'_n$ taken coarsely to some coarse
$\PD(n)$ edge space $X(e)\subset X_n$. The construction of $T'_n$ is a
relative version of the construction in the basis step, which itself is
adapted from the proof of the Homogeneous Theorem~\ref{ThmHomogeneous}.
First one chooses any
$n$-dimensional edge $e$ of $T_n$, between vertices $v',w' \in T'_{n+1}$.
Then for each $h \in H$ one attaches an edge $h \cdot e$ between $h \cdot
v'$ and $h \cdot w'$. The result is not a forest, but it is
quasi-isometrically identified with the forest $T_n$. Now attach 2-cells
in an $H$-equivariant manner, with finitely many $H$-orbits, so that each
of the edges of a given 2-cell have edge spaces in the same coarse
equivalent class. Now we have an $H$-complex each of whose components is
simply connected. Apply tracks to get an $H$-forest containing $T'_{n+1}$
as a subforest. Each new edge (resp.\ vertex) of this forest has an
$n$-dimensional edge space (resp.\ vertex space) mapping back to a
$n$-dimensional edge space of $X_n$.

\section*{Pattern rigidity}

Suppose we are in the setting of Theorem \ref{ThmInhom}. Let
$v'$ be a point raft of $T'$ and $v$ the corresponding point raft of $T$.
Note that as a consequence of Theorem \ref{ThmInhom}, the
quasi-isometry $X'(v') \to X(v)$ takes the edge space pattern inside
$X'(v')$ coarsely to the edge space pattern inside $X(v)$. This information
can be used to strengthen applications of Theorem
\ref{ThmInhom}, by applying ``pattern rigidity'' results. 

For example, in the setting of Theorem \ref{ThmHVZE} we have the
following theorem of R. Schwartz \cite{Schwartz:Symmetric}:

\begin{theorem}[Geodesic Pattern Rigidity in $\hyp^n$] 
\label{ThmGeoPatterns}
Let $G$ be a discrete, cocompact group of isometries of $\hyp^n$, $n \ge
3$. Let $A$ be a nonempty, $G$-equivariant set of geodesics in $\hyp^n$,
with finitely many $G$-orbits. Let $H$ be a group and let $H \cross \hyp^n
\xrightarrow{\phi} \hyp^n$ be a cobounded, proper quasi-action which
coarsely respects $A$: for each $h \in H$ and $a
\in A$ there exists $a' \in A$ such that $\phi(h,a) \ceq a'$. Then there
is an isometric action $\psi \from H \to \Isom(\hyp^n)$ with finite kernel
which strictly respects $A$, such that $\psi$ is a bounded distance from
$\phi$, i.e.\ $\sup_{h,x} d(\phi(h,x),\psi(h,x)) <
\infinity$.
\end{theorem}

In this setting of this theorem, let $\Isom(\hyp^n,A)$ be the set of
isometries of $\hyp^n$ respecting $A$. This is a discrete, cocompact group
of isometries, containing $G$ and $H/\kernel(\Psi)$ as finite index
subgroups. It follows that $G$ and $H/\kernel(\psi)$ are commensurable, by
a commensuration taking the $A$-stabilizers in $G$ to the $A$-stabilizers
of
$H/\kernel(\psi)$. Combined with the discussion above, this completes the
proof of Theorem \ref{ThmHVZE}.

We can get even stronger conclusions under stronger hypotheses. For
instance, consider a graph of groups $\Gamma$ as in Theorem
\ref{ThmHVZE}. Recall that in a discrete, cocompact group of
isometries of $\hyp^n$, the set of loxodromic axis stabilizers is
identical to the set of maximal, virtually cyclic subgroups. Suppose that
we make the following additional assumption:
\begin{itemize}
\item Each edge-to-vertex injection $\xi_\eta \from \Gamma_e
\to\Gamma_{v(\eta)}$ has image equal to a maximal virtually cyclic subgoup
of $\Gamma_{v(\eta)}$, and two distinct edge-to-vertex injections into
$\Gamma_{v(\eta)}$ have distinct images.
\end{itemize}
This implies that for any vertex $v$ of the Bass-Serre tree $T$, distinct
incident edge spaces inside $X(v)$ are all coarsely inequivalent in
$X(v)$. The Bestvina-Feighn Combination Theorem
\cite{BestvinaFeighn:combination} implies that $\pi_1\Gamma$ is word
hyperbolic.

With this additional assumption, Theorem \ref{ThmInhom}
implies that for any group $H$ quasi-isometric to $\pi_1\Gamma$, the
proper, cobounded quasi-action of $H$ on $X$ satisfies tree rigidity.
Combining this with Geodesic Pattern Rigidity, it follows that $H$ is
weakly commensurable to $\pi_1\Gamma$. 

In fact, what this argument shows, under the additional assumption, is that
the group $\pi_1\Gamma$ has finite index in its quasi-isometry group
$\QI(\pi_1\Gamma)$; this is the strongest form of quasi-isometric rigidity.
We also obtain a computation the abstract commensurator of
$\pi_1\Gamma$: it is isomorphic to $\QI(\pi_1\Gamma)$, which is isomorphic
to the full isometry group of the tree of spaces $X$.

\bigskip

Finally we turn to abelian pattern rigidity and a strengthening of
Theorem~\ref{ThmAbelian}. 

\begin{theorem}[Abelian Pattern Rigidity] 
\label{ThmAbelPatterns}
Suppose that $V_1,\ldots,V_K
\subset \Euc^n$ and $W_1,\ldots,W_K \subset \Euc^n$ are affine foliations.
Suppose that there exists a quasi-isometry $f \from \Euc^n \to \Euc^n$
which maps $W_k$ coarsely to $V_k$, for each $k=1,\ldots,K$.
Then there exists a linear isomorphism $F \from \Euc^n \to \Euc^n$ such
that $F(W_k)=V_k$ for each~$k$.
\end{theorem}

\begin{proof}
By passing to asymptotic cones we replace $f$ with a bilipschitz
homeomorphism taking $V_k$ to $W_k$ for each $k$. Applying the Rademacher
Theorem, at almost any point $x$ the derivative
$F=D_xf$ gives the desired conclusion. 
\end{proof}

If $S(V)$ is the linear subspace of $\Euc^n$ parallel to the leaves of an
affine foliation $V$, then a pattern of affine foliations $V_1,\ldots,V_K$
induces a pattern of linear subspaces $S(V_1),\ldots,S(V_K) \subset
\Euc^n$, which in turn induces a pattern of projective subspaces $\P
S(V_1),\ldots,\P S(V_K) \subset \P^{n-1}$, called the \emph{projective
pattern} associated to $V_1,\ldots,V_K$. The Abelian Pattern Rigidity
Theorem \ref{ThmAbelPatterns} shows that the associated projective
pattern is a quasi-isometry invariant of patterns of affine foliations in
$\Euc^n$. For example, four distinct \nb{1}dimensional affine foliations of
$\Euc^2$ have an associated projective pattern of four distinct points in
$\P^1$. The moduli space of such projective patterns is
\nb{1}dimensional, parameterized by the cross-ratio, and so the cross-ratio
is a pattern preserving quasi-isometry invariant.

We can strengthen Theorem \ref{ThmAbelian} as follows. Suppose that
$\Gamma$ is a reduced graph of abelian groups as in Theorem
\ref{ThmAbelian}, with Bass-Serre tree of spaces $X \to T$. Recall that the
hypotheses of Theorem \ref{ThmAbelian} imply that each vertex of $T$ is a
raft. Our strengthened conclusions say that projective patterns associated
to edge space patterns inside vertex spaces of $X$ are
quasi-isometricially rigid. Here are the details.

Let $H$ be a finitely generated group quasi-isometric to $\pi_1\Gamma$.
The proof of Theorem \ref{ThmAbelian} gives a graph of groups $\Gamma'$
with Bass-Serre tree of spaces $X' \to T'$, such that $H
\approx\pi_1\Gamma'$, together with a quasi-isometry $\phi \from X' \to
X$, such that for each edge $e'$ of $T'$ there is an edge $e$ of $T$ with
$\phi(X(e)) \ceq X'(e')$, and for each vertex $v'$ of $T'$, either $v'$ is
a raft of $T'$ and $\phi(X'(v'))\ceq X(v)$ for some vertex $v$ of $T$, or
$v'$ is not a raft of $T'$ and so $v'$ has an incident edge $e'$ with
$X'(v') \ceq X'(e')$. 

Consider a raft vertex $v'$ of $T'$ and the corresponding vertex $v$ of
$T$. Composing the map $X(v') \xrightarrow{\phi} X'$ with the closest point
projection to $X(v)$ gives a quasi-isometry still denoted $\phi \from
X'(v') \to X(v)$. Moreover, this quasi-isometry takes the edge space
pattern inside $X'(v') \approx \Euc^n$ coarsely to the edge space pattern
inside $X(v) \approx \Euc^n$. Applying the Abelian Pattern Rigidity
Theorem \ref{ThmAbelPatterns}, the projective patterns in
$\P^{n-1}$ associated to $X'(v')$ and $X(v)$ are projectively equivalent.


\providecommand{\bysame}{\leavevmode\hbox to3em{\hrulefill}\thinspace}

\bigskip

{\small
\noindent
Lee Mosher:\\
Dept.\ of Mathematics and Computer Science\\
Rutgers University, Newark\\
Newark, NJ 07102\\
E-mail: mosher@andromeda.rutgers.edu

\bigskip\noindent
Michah Sageev:\\
Dept. of Mathematics\\
Technion---Israel Institute of Technology\\
Haifa 32000, Israel\\
E-mail: sageevm@techunix.technion.ac.il

\bigskip\noindent
Kevin Whyte:\\
Dept.\ of Mathematics\\
Univ. of Utah\\
Salt Lake City, Utah 84112-0090\\
E-mail: kwhyte@math.utah.edu
}

\end{document}